\documentclass[number,seceqn,dvips]{arxbj}
\usepackage{mathbh}

\aid{0}
\volume{16}
\issue{4}
\pubyear{2010}
\firstpage{909}
\lastpage{925}
\doi{10.3150/09-BEJ236}

\makeatletter

\newtheorem{theorem}{Theorem}[section]
\newtheorem{itlemma}[theorem]{Lemma}
\newtheorem{itproposition}[theorem]{Proposition}

\newremark{itremark}[theorem]{Remark}

\makeatother

\begin{document}
\begin{frontmatter}

\title{Tightness for the interface of the one-dimensional contact process}
\runtitle{Tightness for the interface of the contact process}

\begin{aug}
\author[a]{\fnms{Enrique} \snm{Andjel}\thanksref{a}\ead[label=e1]{enrique.andjel@cmi.univ-mrs.fr}},
\author[b]{\fnms{Thomas} \snm{Mountford}\thanksref{b}\ead[label=e2]{thomas.mountford@epfl.ch}},

\author[c]{\fnms{Leandro P.R.} \snm{Pimentel}\thanksref{c}\ead[label=e3]{leandro@im.ufrj.br}}
\and
\author[b]{\fnms{Daniel} \snm{Valesin}\thanksref{b}\ead[label=e4]{daniel.valesin@epfl.ch}\corref{}}
\runauthor{Andjel, Mountford, Pimentel and Valesin}
\address[a]{Centre de Math\'ematiques et Informatique (CMI),
Universit\'e de Provence,
36 Rue Joliot-Curie,\\ 13013~Marseille, France. \printead{e1}}
\address[b]{Institut de Math\'ematiques, Station 8,
\'Ecole Polytechnique F\'ed\'erale de Lausanne,
CH-1015 Lausanne, Switzerland. E-mails: \printead*{e2,e4}}
\address[c]{Institute of Mathematics,
Federal University of Rio de Janeiro,
IM-UFRJ, Postal Code 68530,\\ CEP 21945-970, Rio de Janeiro, RJ, Brazil.
\printead{e3}}
\end{aug}

\received{\smonth{3} \syear{2009}}
\revised{\smonth{9} \syear{2009}}

%
\begin{abstract}
We consider a symmetric, finite-range contact process with two types of
infection; both have the same (supercritical)
infection rate and heal at rate 1, but sites infected by Infection 1
are immune to Infection 2. We take the initial
configuration where sites in $(-\infty, 0]$ have Infection 1 and sites
in $[1, \infty)$ have Infection 2, then
consider the process $\rho_t$ defined as the size of the interface
area between the two infections at time~$t$.
We show that the distribution of $\rho_t$ is tight, thus proving a
conjecture posed by Cox and Durrett in
[\textit{Bernoulli}~\textbf{1} (1995) 343--370].
\end{abstract}

%

\end{frontmatter}

\section{Introduction}
\label{S1}

This paper addresses a conjecture of Cox and Durrett \cite{cd}
concerning interfaces naturally arising in supercritical contact
processes on $\mathbb{Z}^1$.

The contact process on $\mathbb{Z}^d$ is a spin system with operator
\[
\Omega f(\eta) = \sum_x \bigl(f (\eta^{x}) - f (\eta
)\bigr) c(x, \eta
), \qquad \eta\in{\{0, 1\}}^{\mathbb{Z}^d},
\]
where
\[
\cases{ \eta^x (y) = \eta(y) , &\quad\mbox{if } $y \neq x$,
\cr
\eta^x(x) = 1-\eta(x),&
}
\]
and flip rates $c(x, \eta)$ are given by
\[
c(x, \eta) = \cases{
 1 , &\quad\mbox{if } $\eta(x) = 1 ,$
 \cr
\lambda\sum p(y-x) \eta(y) , &\quad\mbox{if } $\eta(x) = 0 $
}
\]
for $\lambda> 0$ and probability kernel $p(\cdot)$.

In the following, we take $p(\cdot)$ to have finite range (that is,
$\exists M < \infty  \dvtx p(x) = 0$ for $|x| > M$) and to be symmetric,
though this latter hypothesis can be dispensed with via the techniques
and results of Bezuidenhout and Gray \cite{bgray}.

Often the contact process is used as a model of the spread of an
infection and a configuration $\eta\in{\{0,1\}}^{\mathbb{Z}^{d}}$ represents
the state where there is an infection at $x \in\mathbb{Z}^{d}$ if and
only if
$\eta(x) = 1$. We will adopt this point of view and speak of a site
$x$ being \textit{infected} at time $t$ (for a process $(\eta_{t} \dvtx t
\geq0))$ if $\eta_{t}(x) = 1$. We will sometimes identify
configurations in $\{0,1\}^{\mathbb{Z}^d}$ with their sets of infected sites
(that is, we will write $\xi$ instead of $\{x\dvtx  \xi(x) = 1\})$. As
defined above, the contact process is attractive (see \cite{lig85} for
fundamental results associated with this property). Thus, for two
configurations $\xi_{0} $ and $\zeta_{0}$ satisfying $\xi_{0} \leq
\zeta_{0}$ under the natural partial order, it is possible to
construct in a single probability space two processes, $(\xi_{t} \dvtx t
\geq0)$ starting at $\xi_{0}$ and $(\zeta_{t} \dvtx t \geq0)$ starting
at $\zeta_{0}$, satisfying, with probability one, $\xi_{t} \leq\zeta
_{t}$ for all $t$.

A consequence is that $\exists  \lambda^1_c$ such that for $\lambda
> \lambda^1_c$, the invariant limit $\lim_{t \rightarrow\infty}
\delta_{\mathbh{1}}S(t)$ is a non-trivial measure and for $\lambda<
\lambda_{c}^{1}$, this limit is $\delta_{\underline{0}}$. There also
exists $\lambda_{c}^{2}$ such that for $\lambda> \lambda_{c}^{2},
P^{\{0\}} (\tau= \infty) > 0$ for $\tau= \inf\{t\dvtx  \eta_{t} \equiv
0\}$, and for $\lambda< \lambda_{c}^{2},  P^{\{0\}}(\tau= \infty)
= 0$. In fact, via duality (see, for example, \cite{d1} or \cite{lig85}), $\lambda^{1}_{c} = \lambda^{2}_{c}$, and this critical value
will henceforth be denoted by $\lambda_{c}$.

We now introduce some notation. Suppose we are given independent
Poisson processes on $[0, \infty)$, $\{D_{x}\}_{x \in\mathbb{Z}^{d}}
$ of
rate 1 and $\{N^{(x, y)}\}_{x, y \in\mathbb{Z}^{d}}$ of rate $\lambda
p(y-x)$. Denote by $H$ a realization of all these independent
processes; we say that $H$ is a \textit{Harris construction}. $H$ is
thus a Poisson measure on $(\mathbb{Z}^d \cup(\mathbb{Z}^d)^2)
\times[0, \infty)$
such that, if $y, z \in\mathbb{Z}^d$ and $I$ is a Borelian subset of $[0,
\infty)$, we have $H(\{z\} \times I) = D_z(I)$ and $H(\{(y, z)\}\times
I) = N^{(y, z)}(I)$. Given a Harris construction $H$ and $(x, t) \in
\mathbb{Z}
^d \times[0, \infty)$, denote by $H^{(x, t)}$ the Harris construction
obtained by shifting $H$ so that the space origin becomes $x$ and the
time origin becomes $t$. Formally, if $y, z \in\mathbb{Z}^d$ and $I$
is a
Borelian subset of $[0, \infty)$, then $H^{(x, t)}(\{z\}\times I) =
H(\{z+x\}\times(I+t))$ and $H^{(x, t)}(\{(y, z)\} \times I) = H(\{
(y+x, z+x)\}\times(I + t))$.

Given a Harris construction $H = \{(D_x)_{x \in\mathbb{Z}^d}, (N^{(x,y)})_{x,
y \in\mathbb{Z}^d}\}$ and $(x, s), (y, t) \in\mathbb{Z}^d \times
\mathbb{R}_+$ with $s <
t$, we write $(x, s) \leftrightarrow(y, t)$ (in $H$) if there exists a
piecewise constant $\gamma\dvtx [s, t] \rightarrow \mathbb{Z}^d$ such
that:
\begin{longlist}[(iii)]
\item[(i)] $\gamma(s) = x,  \gamma(t) = y;$

\item[(ii)] $\gamma(r) \neq\gamma(r-) \mbox{ only if } r \in N^{\gamma
(r-), \gamma(r)};$

\item[(iii)] $\not\exists  s \leq r \leq t \mbox{ with } r \in D_{\gamma
(r)}.$
\end{longlist}
Given $A, B, C \subset\mathbb{Z}^d$ and $s, t \in\mathbb{R}_+$, we
write $A \times
s \leftrightarrow B \times t$ if $(x, s) \leftrightarrow(y, t)$ for
some $x \in A$, $y \in B$. Additionally, $A \times\{s\} \leftrightarrow
B \times\{t\}$ \emph{inside} $C$ if there exists a path connecting $A
\times\{s\}$ and $B \times\{t\}$ and with image contained in $C$.

Given $\xi_0 \in\{0, 1\}^{\mathbb{Z}^d}$ and a Harris construction
$H$, we
construct a trajectory $(\eta_t^{\xi_0}(H) \dvtx t \geq0)$ by specifying
$\eta_0^{\xi_0}(H) = \xi_0$ and $[\eta_{t}^{\xi_0}(H)](x) = 1$ if
and only if $\xi_0 \times\{0\} \leftrightarrow(x, t)$ in $H$.

A moment's reflection shows that, under the law of $H, (\eta_t^{\xi
_0}(H))_{t \geq0}$ is a contact process with initial condition $\xi
_0$ and, if $\xi_0 \leq\zeta_0$, then putting $\xi_t = \eta_t^{\xi
_0}(H)$ and $\zeta_t = \eta_t^{\zeta_0}(H)$, we obtain the claimed
coupling of two processes, one of which is always inferior to the other.

As noted, we will be concerned with one-dimensional contact processes
with $\lambda> \lambda_{c}$. Define $r_t^{\xi_0}(H) = \sup\{x\dvtx
[\eta_t^{\xi_0}(H)](x) = 1\}$. We will usually omit the dependency on
$H$ and when we omit the initial condition and simply write $r_t$, we
take $\xi_0 = I_{(-\infty, 0]}$. If $\xi_0$ is such that $\sum_x
\xi_0(x) = \infty$ and $\sup\{x\dvtx  \xi_0(x)=1\} < \infty$, then
almost surely $\eta_t^{\xi_0} \neq\underline{0}$ and $r_t^{\xi_0}
< \infty$ for all $t$. It is classical that $\frac{r_t}{t} \stackrel
{t \to\infty}{\longrightarrow} \alpha= \alpha(\lambda) > 0$; see
Theorems 2.19 and 2.27 in \cite{lig85} (even though the process
treated there is nearest-neighbor, the proof works for the finite-range
case as well).

We consider the following question. Define
\begin{eqnarray*}
l_t &=& l_t(H) = \inf\bigl\{x\dvtx  \bigl[\eta_t^{(-\infty, 0]}(H)\bigr](x) \neq[\eta
_t^{\mathbh{1}}(H)](x)\bigr\},
\\
\rho_t &=& r_t^{(-\infty, 0]} - l_t,\qquad  \rho_t^+ = \max\{\rho_t, 0\}
, \qquad \rho_t^-=\max\{-\rho_t,0\}.
\end{eqnarray*}
While it is easy to see that $\{r_t < l_t\}$ and $\{l_t < r_t\}$ are
events of strictly positive probability, it is reasonable to believe
that the two quantities are close. Cox and Durrett conjectured that ${\{
\vert\rho_t \vert\}}_{t \geq0}$ would be a tight collection of
random variables. We answer the conjecture affirmatively.

\begin{theorem}\label{tightness}
The law of $\{\rho_t\}_{t\geq0}$ is tight. That is, for any $\delta>
0$, there exists $L > 0$ such that $\mathbb{P}(|\rho_t|>L)< \delta$
for every $t \geq0$.
\end{theorem}

From the joint process $((\eta_{t}^{\mathbh{1}}, \eta_{t}^{(-\infty,
0]}) \dvtx  t \geq0)$, we can define a process $(\chi_{t} \dvtx t \geq0)$ on
$ {\{0, 1, 2\}}^{\mathbb{Z}}$ by
\[
\chi_{t} (x) =\cases{ 0 , &\quad\mbox{if } $\eta_{t}^{(-\infty, 0]}(x) =
\eta
_{t}^{\mathbh{1}}(x) = 0 $,
\cr
1 , &\quad\mbox{if } $\eta_{t}^{(-\infty, 0]}(x) = 1 $,
\cr
2,  &\quad\mbox{if } $\eta_{t}^{(-\infty, 0]}(x) = 0,  \eta_{t}^{\mathbh{1}}(x) = 1.$
}
\]
It is not difficult to see that $\chi_{t}$ is a realization of a
process taking values in $\{0, 1, 2\}^\mathbb{Z}$ with initial configuration
equal to $I_{(-\infty, 0]} + 2\cdot I_{(0, \infty)}$ and the
following rates:
\begin{eqnarray*}
0 &\to&1\qquad \mbox{at rate } \lambda\sum p(y-x) I_{\chi(y) = 1};\\
0 &\to&2\qquad \mbox{at rate } \lambda\sum p(y-x) I_{\chi(y) = 2};\\
2 &\to&0\qquad \mbox{at rate } 1;\\
1 &\to&0\qquad \mbox{at rate } 1;\\
2 &\to&1\qquad \mbox{at rate } \lambda\sum p(y-x) I_{\chi(y) = 1}.
\end{eqnarray*}

The particle system with the above transition rates is a model for
hierarchical competition considered in \cite{durmol} and \cite
{durswin}; the following interpretation is provided. Sites in state 0
are said to contain grass, in state 1 to contain trees and in state 2
to contain bushes. When trees attempt to occupy new territory, they are
able to displace bushes, but bushes cannot displace trees. Since, in
our case, we take the initial configuration $I_{(-\infty, 0]} + 2\cdot
I_{(0, \infty)}$, we expect the area taken by trees to grow to the
right towards the area originally taken by bushes. However, since we
allow for non-nearest-neighbor interactions, we may observe a mixed
area where the two coexist. With the above notation, this area appears
when $\rho_t > 0$. Alternatively, it may happen that there is no mixed
area and a gap of grass appears between the two homogeneous zones (in
the case $\rho_t < 0$). Theorem \ref{tightness} states that with
large probability, and uniformly in time, neither the mixed nor the
intermediate grass area is too large.

The proof is divided into two parts. The first part, namely the proof
of tightness of $\{\rho_t^+\}$, is given at the end of Section \ref{S2}. The
key ingredients are the celebrated result of Bezuidenhout and Grimmett \cite{bg},
the renormalization arguments employed by, among others, Durrett (see
\cite{d1}) and the construction carried out in \cite{MS}. These
permit us to argue that from a single $(x, t)$ with $\eta
_{t}^{(-\infty, 0]}(x) = 1$, there will be positive probability that
inside a cone $C_{x, t} = \{(y, s) \dvtx \vert y -x \vert\leq\beta(s-t)\}
$, $\eta^{\mathbh{1}}$ and $\eta^{(-\infty, 0]}$ are equal. In
Section 3, a much simpler argument is employed to establish tightness
of $\{\rho_t^-\}$.

\section{Tightness of $\{\rho_t^+\}$}
\label{S2}

\subsection{Right edge speed}
\label{S2.1}

Given $\gamma> 0$, we say that $(0, 0)\in\mathbb{Z}\times[0, \infty
)$ is
$\gamma$\emph{-slow up to time} $T$ if $r_t \leq\gamma t \,  \forall
t \leq T$. If this is satisfied for all $T$, then we say that $(0, 0)$
is $\gamma$\emph{-slow}.

\begin{itlemma}\label{gamma}
\textup{(i)} For any $\varepsilon> 0$, there exists $\gamma> 0$
such that $\mathbb{P}((0,0) \mbox{ is } \gamma\mbox{-slow}) > 1 -
\varepsilon$.\vspace*{-6pt}
\begin{longlist}
\item[(ii)] For any $\gamma> \alpha$, we have
%
\begin{equation}\label{pslow} \mathbb{P}((0,0) \mbox{ is } \gamma
\mbox{-slow}) > 0
\end{equation}
and there exist $c, C > 0$ such that
%
\begin{equation}\label{expslow} \mathbb{P}((0, 0) \mbox{ is } \gamma
\mbox{-slow up to time }
T \mbox{ but not } \gamma\mbox{-slow}) \leq C \mathrm{e}^{-cT}.
\end{equation}
\end{longlist}
\end{itlemma}

\begin{pf} Almost surely, $t \mapsto r_t$ is right-continuous with left
limits, identically zero in a neighborhood of 0 and satisfies $r_t/t
\rightarrow \alpha$. It follows that almost surely, $\{r_t/t\dvtx t \geq
0\}$ is
bounded, hence we have (i). It also follows that, given $\gamma>
\alpha$, we can obtain $R > 0$ such that \mbox{$\mathbb{P}(r_t/t < R/t +
\gamma\,
\forall t) > 0$}. Now,
\[
\mathbb{P}((0,0) \mbox{ is } \gamma\mbox{-slow}) \geq\mathbb
{P}\bigl(r_t \leq0
\,\forall t \in[0, R/\gamma], r_s^{(-\infty,0]}\bigl(H^{(r_{R/\gamma},
R/\gamma)}\bigr) < R + \gamma s  \, \forall s \geq0\bigr).
\]
The first event on the above probability depends only on the Harris
construction $H$ on $[0, R/\gamma]$, whereas the second depends only
on $H$ on $[R/\gamma, +\infty)$, so they are independent. Also noting
that $\mathbb{P}(r_s^{(-\infty,0]}(H^{(r_{R/\gamma}, R/\gamma)}) <
R +
\gamma s  \, \forall s \geq0) = \mathbb{P}(r_s < R + \gamma s
\,\forall s
\geq0),$ we get, by translation invariance,
\[
\mathbb{P}((0,0) \mbox{ is } \gamma\mbox{-slow}) \geq\mathbb
{P}(r_t \leq0
\,\forall t \in[0, R/\gamma]) \cdot\mathbb{P}(r_s < R + \gamma s
\,\forall s
\geq0).
\]
The second probability above is positive by our choice of $R$. The
first one is also positive because it contains the event $\{(-\infty,
0] \times[0, R/\gamma] \nleftrightarrow(0, +\infty) \times[0,
R/\gamma]\}$, which has positive probability since it corresponds to a
finite number of Poisson processes having no arrivals in a finite time
interval. We thus have (\ref{pslow}).

To establish (\ref{expslow}), fix $\gamma' \in(\alpha, \gamma)$
and note that
\begin{eqnarray*}
&&\mathbb{P}(r_t \leq\gamma t \mbox{ for all } t \in[0, T] \mbox{
but not
for all } t \geq0)
\\
&&\quad \leq\mathbb{P}(\exists t > T\dvtx r_t > \gamma t)
\leq\mathbb{P}(\exists t > T\dvtx r_t > \gamma t, r_T \leq\gamma' T) +
\mathbb{P}(r_T
> \gamma'T) .
\end{eqnarray*}

By Lemma 2 in \cite{MS} (a large deviations result for $r_t$), $\gamma
' > \alpha$ implies that the second term in the sum decays
exponentially fast in $T$ and, by translation invariance, the first
term is less than $\mathbb{P}(\exists s> 0: r_s > (\gamma- \gamma')
T +
\gamma s)$. It will therefore suffice to prove that $\mathbb
{P}(\exists s > 0\dvtx
r_s > k + \gamma s)$ decays exponentially fast as $k$ tends to
infinity. Indeed, put $\theta= \mathbb{P}(\exists t > 0\dvtx r_t \geq M +
\gamma
t)$ (remember that $M$ is the range of the process) and $T_N = \inf\{t
\geq0: r_t \geq2MN + \gamma t\}$ for $N \geq1$. We have $\theta< 1$
by (\ref{pslow}) and
\begin{eqnarray*}
\mathbb{P}(T_{N+1}<\infty) &=&\mathbb{P}\bigl(\exists t > 0\dvtx r_t \geq
2M(N+1) + \gamma t\bigr)
\\
&\leq&\mathbb{P}\bigl(T_N < \infty, \exists s > 0\dvtx r_s^{(-\infty, 0]}\bigl(H^{(r_{T_N},
T_N)}\bigr) \geq M + \gamma s\bigr)
\\
&=&\mathbb{P}(T_N < \infty) \cdot\mathbb{P}(\exists s > 0\dvtx r_s > M +
\gamma s) = \mathbb{P}
(T_N < \infty) \cdot\theta.
\end{eqnarray*}

Thus $\mathbb{P}(T_N < \infty) \leq\theta^N$. Now, if $k \geq1$, then
\[
\mathbb{P}(\exists s > 0\dvtx r_s > k + \gamma s) \leq\mathbb{P}(\exists
s> 0\dvtx r_s \geq
2M\sigma+ \gamma s) \leq\mathbb{P}(T_\sigma< \infty) \leq\theta
^\sigma,
\]
where $\sigma$ denotes the largest integer strictly smaller than $k/2M$.
\end{pf}

\subsection{Descendancy barriers}
In this section, we define an event called the \emph{formation of a
descendancy barrier}. This will mean that, inside a certain area
delimited by a vertical cone that grows upward from the origin, all
infected sites will be connected to the origin. Additionally, no
infection from one side of the cone will be able to pass to the other
side without being connected to the origin. These barriers, which
appear with positive probability, as we will show, are the essential
structure in our proof of tightness of $\{\rho_t^+\}$.

We first give a brief exposition of oriented percolation and state a
result that will be needed later. For a detailed treatment of the
subject, see the survey \cite{d2}.

Let $\Lambda= \{(m, n) \in\mathbb{Z}\times\mathbb{Z}_+\dvtx m+n \mbox
{ is even}\},
\Omega= \{0,1\}^\Lambda$ and $\mathcal{F}$ be the $\sigma$-algebra
generated by cylinder sets of $\Omega$. Points of $\Omega$ will be
denoted by $\Psi$, with $\Psi(m,n) \in\{0, 1\}$ for $(m, n) \in
\Lambda$. $\mathbb{P}_p$ will denote the product measure $(p \delta
_1 +
(1-p)\delta_0)^{\otimes\Lambda}$. The vertical axis of $\Lambda$
will be interpreted as time.

Given $k \geq1,  \varepsilon> 0$ and a probability $\mathbb{P}$ on
$\mathcal{F}$, we
say that $(\Omega, \mathcal{F}, \mathbb{P})$ is a $k$-dependent
oriented percolation
system with closure below $\varepsilon$ if

\begin{equation}\label{kdepper}
\mathbb{P}\bigl(\Psi(m_i, n) = 0, 1 \leq i \leq r \mid\{\Psi(m, s)\dvtx 1
\leq s <
n, (m, s) \in\Lambda\}\bigr) < \varepsilon^r ,
\end{equation}
where $r \geq1$, $(m_i, n) \in\Lambda \, \forall i$ and $|m_{i_1} -
m_{i_2}| > 2k$ when $i_1 \neq i_2$.

Given $\Psi\in\Omega$, we say that two points $(x,m), (y,n) \in
\Lambda$ with $m < n$ are \emph{connected by an open path} if there
exists a sequence $x_0 = x, x_1, \ldots, x_{n-m} = y$ in $\mathbb{Z}$ such
that $|x_{i+1}-x_i| = 1 \,  \forall i \in\{0, \ldots, n-m-1\}$ and
$\Psi(x_i, m+i) = 1 \,  \forall i \in\{0, \ldots, n-m\}$. We say
that $(x, m)$ \emph{percolates up to time} $n$ when it is connected by
an open path to a point at height $n$. Finally, we say that $(x, m)$
\emph{percolates} when there is an infinite open path starting from it.

In \cite{d2}, it is proved that if $p$ is sufficiently large, then the
origin percolates with positive probability in $(\Omega, \mathcal{F},
\mathbb{P}
_p)$. Moreover, the rightmost particle connected to the origin at time
$n$, denoted $R_n$, almost surely satisfies $\lim R_n/n = \tilde\alpha
(p) > 0$ as $n \rightarrow\infty$. To obtain similar results for
$k$-dependent systems, we use the following particular case of Theorem
0.0 in \cite{lss}.

\begin{itlemma}\label{lss}
Fix $k \in\mathbb{N}$ and $0 < p < 1$. There exists $\varepsilon> 0$
such that
if $(\Omega, \mathcal{F}, \mathbb{P})$ is a $k$-dependent oriented
percolation
system with closure below $\varepsilon$, then $\mathbb{P}$ stochastically
dominates $\mathbb{P}_p$.
\end{itlemma}

Using these facts and an argument similar to the one used in Lemma \ref
{gamma}, we can prove
the following lemma.

\begin{itlemma}\label{betap} Fix an arbitrary $0 < \beta< 1$ and
define the events
\begin{eqnarray}\label{percevents}
\Gamma(i) &=&\left\{ \matrix{\mbox{There exist two open paths, one
starting at } (-2,0) \mbox{, the other at } (2, 0)
\cr  \mbox{ and both
reaching time } i. \mbox{ Neither of them intersects }\{(m, n)\dvtx -\beta
n \leq m \leq\beta n \}
}\right\},\qquad\nonumber
\\[-8pt]\\[-8pt]
\Gamma &=& \left\{\matrix{
\mbox{There exist two infinite open paths,
one starting at } (-2,0) \cr
\mbox{and the other at } (2, 0). \mbox{ Neither of them intersects }\{(m, n): -\beta n \leq m \leq\beta n\}
}\right\}.\nonumber
\end{eqnarray}
For any $k$ and $\bar\delta>0$, there exists $\varepsilon>0$ such that
if $(\Omega, \mathcal{F}, \mathbb{P})$ is a $k$-dependent
percolation system with
closure below $\varepsilon$, then:
\begin{longlist}[(ii)]
\item[(i)] $\mathbb{P}(\Gamma) > 1-\bar\delta$;
\item[(ii)] $\mathbb{P}(\Gamma(i) \backslash\Gamma) \leq D
\mathrm{e}^{-d i}$
for some $d, D > 0$.
\end{longlist}
\end{itlemma}

We now construct a mapping $H \mapsto\Psi_H$ of Harris constructions
into points of $\Omega$; this is essentially a repetition of the
mapping developed in \cite{MS}. The construction will depend on large
integers $K$ and $N$ (in particular, much larger than the range $M$)
whose choice will be described in Proposition \ref{MSprop}. Given $m
\in\mathbb{Z}, n \in\mathbb{Z}_+$, define
\begin{eqnarray}\label{ImJmn}
I_m &=& \biggl(\frac{mN}{2} - \frac{N}{2}, \frac{mN}{2} + \frac
{N}{2}\biggr] \cap\mathbb{Z} ,\nonumber
\\[-8pt]\\[-8pt]
J_{(m, n)} &=& \biggl[\frac{mN}{2} - M, \frac{mN}{2} + M \biggr]
\times[KNn, KN(n+1)] \cap\mathbb{Z}\times[0, +\infty) .\nonumber
\end{eqnarray}
We start defining an auxiliary $\Phi_H \in\{0,1,2\}^\Lambda$. Given
$(m, 0) \in\Lambda$, put $\Phi_H(m, 0) = 1$ if $H$ and the
trajectory $\eta^{\mathbh{1}}(H)$ satisfy the following conditions:
%

\begin{eqnarray}\label{0mscond1}
\mbox{there is no vacant interval at time } KN
\mbox{ of length }
N^{1/2} \mbox{ inside } I_{m-1}\cup I_{m+1};
\end{eqnarray}
%
\begin{eqnarray}\label{0mscond2}
\mbox{every occupied site in }I_{m-1} \cup I_{m+1} \mbox{ at time }KN
\mbox{ is a descendant of }I_m \times\{0\};
\end{eqnarray}
\begin{eqnarray}\label{0mscond3}
&&\mbox{there does not exist }(z, s) \in J_{(m, 0)}\mbox{ such that }\nonumber
\\[-11pt]\\[-11pt]
 &&I_m \times\{0\} \nleftrightarrow(z, s) \mbox{ and } (I_m^C \times
[0, s]) \leftrightarrow(z, s);\nonumber
\end{eqnarray}
put $\Phi_H(m,0) = 0$ otherwise. Given $(m, n) \in\Lambda$ with $n
\geq1$, put $\Phi_H(m,n) = 1$ if
%
%
\begin{eqnarray}\label{mscond1}
1 \in\{\Phi_H(m-1,n-1), \Phi_H(m+1,n-1)\};
\end{eqnarray}
%
\begin{eqnarray}\label{mscond11}
&&\mbox{there is no vacant interval at time } KN(n+1)\nonumber
\\[-11pt]\\[-11pt]
&&\mbox{of
length } N^{1/2} \mbox{ inside } I_{m-1}\cup I_{m+1};\nonumber
\end{eqnarray}
\begin{eqnarray}\label{mscond2}
&&\mbox{every occupied site in }I_{m-1} \cup I_{m+1} \mbox{ at time
}KN(n+1)\nonumber
 \\[-11pt]\\[-11pt]
&&\mbox{is a descendant of }(I_m \cap\eta^\mathbh
{1}_{KNn})\times KNn;\nonumber
\end{eqnarray}
\begin{eqnarray}\label{mscond3}
&&\mbox{there does not exist }(z, s) \in J_{(m, n)}\mbox{ such that }\nonumber
\\[-11pt]\\[-11pt]
 &&\bigl((I_m \cap\eta^{\mathbh{1}}_{KNn}) \times KNn\bigr) \nleftrightarrow(z,
s) \mbox{ and } (I_m^C \times[KNn, s]) \leftrightarrow(z, s).\nonumber
\end{eqnarray}
If (\ref{mscond1}) fails, put $\Phi_H(m, n) = 2$, and in every other
case, put $\Phi_H(m, n) = 0$. Finally, set
\[
\Psi_H(m,n) = \cases{ 0 , &\quad\mbox{if } $\Phi_H(m,n) = 0 $,
\cr
 1 , &\quad\mbox{otherwise.}
}
\]

Note that, with this construction, if there is an infinite open path $\{
(m_i, n_i)\}_{i \geq0}$ leaving the origin in $\Psi_H$, we must have
$\Phi_H(m_i, n_i) = 1$ for every $i$.

We now have the
following proposition.
\begin{itproposition}[(Mountford and Sweet \cite{MS})]\label{MSprop} There exist $k, K$
-- depending only on the parameter $\lambda$ of the contact process --
with the following property: for any $\varepsilon> 0$, there exists $N$
such that $\Psi_H$ defined from $K$ and $N$ is a $k$-dependent
percolation system with closure below $\varepsilon$.
\end{itproposition}

\begin{itremark}\label{C}Conditions (\ref{0mscond1}) and (\ref
{mscond11}) are only necessary
to establish Proposition \ref{MSprop} and will not be used in the
sequel. Also, $N$ in Proposition \ref{MSprop} can be chosen as large
as we want; in particular, as already mentioned, we take both $K$ and
$N$ to be larger than the range $M$.
\end{itremark}

In what follows, the oriented percolation dependency parameter $k$, the
constant $\beta$ and associated events $\Gamma, \Gamma(i)$, the
renormalization constants $N, K$ and the closure density $\varepsilon$
will be fixed in the following way:

\begin{itemize}
\item[$\bullet$]$K$ and $k$ are functions of $\lambda$, as explained
in the last proposition above;
\item[$\bullet$] $\beta$ will be any fixed number in $(0, 1)$;
\item[$\bullet$] $\Gamma$ and $\Gamma(i)$ will be defined from
$\beta$, as in (\ref{percevents});
\item[$\bullet$] $\delta> 0$ will be given during the proof of
Theorem \ref{tightness};
\item[$\bullet$] $\varepsilon$ will be chosen corresponding to $\bar
\delta=\delta/6,  k,  \beta$, as in Lemma \ref{betap};
\item[$\bullet$] $N$ will be chosen corresponding to $\varepsilon$, as
in Proposition \ref{MSprop}.
\end{itemize}

Introducing some more terminology, we call the origin $\beta$-\emph
{expanding} when:
%
\begin{equation}\label{beta1}
\mbox{If } x \in\mathbb{Z} , y \in I_{-2} \cup I_0 \cup I_2 , x
\neq y ,
t \leq1 \mbox{ and } (x,0) \leftrightarrow(y, t) \mbox{, then }
(0,0) \leftrightarrow(y, t);\vspace*{-13pt}
\end{equation}
%
\begin{eqnarray}\label{betanew}
D_0 \cap[0, 1] &=& \varnothing;
\\
\label{beta2}
(0,0) \leftrightarrow(z, 1) \,  \forall z &\in& I_{-2} \cup I_0 \cup I_2;
\\
\label{beta3}
\Psi_{H^{(0,1)}} &\in& \Gamma.
\end{eqnarray}

Condition (\ref{beta1}) means that whenever an infection is
transmitted to a site in $I_{-2} \cup I_0 \cup I_2$ before time $1$,
there must exist an earlier/simultaneous (possibly indirect)
transmission from $(0, 0)$ to the same site. Condition~(\ref{betanew})
means that there is no healing at $\{0\} \times[0,1]$. Condition~(\ref
{beta2}) means that at time $1$, every site in $I_{-2}\cup I_0 \cup
I_2$ carries an infection that descends from the origin.
Condition~(\ref{beta3}) states that the percolation structure defined after
placing the origin at $(0, 1)$ has the properties defined in (\ref
{percevents}). The $\beta$ dependency is in the third event since
$\Gamma$ depends on $\beta$, and also in the choice of the parameters
of the renormalization.

We say that $(0, 0)$ is $\beta$-\emph{expanding up to a time} $T > 1$
when (\ref{beta1})--(\ref{beta2}) are satisfied
and $\Psi_{H^{(0,1)}} \in\Gamma(i)$, where $i$ satisfies $T \in(1 +
KN(i-1), 1 + KNi]$. We then have the
following lemma.
\begin{itlemma}\label{beta}
\textup{(i)} $\mathbb{P}((0, 0) \mbox{ is } \beta\mbox
{-expanding}) > 0$.\vspace*{-6pt}
\begin{longlist}[(ii)]
\item[(ii)] $\mathbb{P}((0, 0) \mbox{ is } \beta\mbox
{-expanding up
to time } T \mbox{, but not } \beta\mbox{-expanding}) \leq\bar D
\mathrm{e}^{-\bar d T}$ for some $\bar d, \bar D > 0$.
\end{longlist}
\end{itlemma}

\begin{pf} It is clear that with positive probability, (\ref
{beta1})--(\ref{beta2}) happen simultaneously. Also, they are
independent of (\ref{beta3}), which, in turn, has positive
probability, by Lemma \ref{betap}, since $\Psi_H$ is supercritical.
Hence, the origin has positive probability of being $\beta$-expanding,
proving (i). Now, note that
\begin{eqnarray*}
&&\{(0,0) \mbox{ is } \beta\mbox{-expanding up to time } T \mbox{,
but not } \beta\mbox{-expanding}\}
\\
&&\quad \subset\bigl\{\Psi_{H^{(0,1)}} \in
\Gamma\bigl(\lfloor(T-1)/KN\rfloor\bigr) \backslash\Gamma\bigr\},
\end{eqnarray*}
where $\lfloor x \rfloor$ denotes the integer part of $x$. The
probability of the last event in the above expression is bounded by $D
\mathrm{e}^{-d(\lfloor(T-1)/KN\rfloor)}$, by Lemma \ref{betap}, so we have (ii).
\end{pf}

Let us now present the properties that motivated this construction. We
start defining, for $\rho>0$,
\[
V(\rho)= \{(z, s) \in\mathbb{Z}\times[0, +\infty)\dvtx - \rho s \leq z
\leq
\rho s\}.
\]
We then have the following proposition.
\begin{itproposition}\label{barrier}
Suppose that the origin is $\beta$-expanding. There then exists a
(deterministic) $0< \bar\beta< 1$ with the following three
properties:
\begin{longlist}[(iii)]
\item[(i)] if $x, z \in\mathbb{Z},  (x, 0) \leftrightarrow(z, s)$ and $(z,
s) \in
V(\bar\beta)$, then $(0,0) \leftrightarrow(z, s)$;
\item[(ii)] $r_s^0 \geq
\cases{ \bar\beta s , \mbox{if } s\geq1 \cr 0 , \mbox
{ if } s<1
}
 \geq\max\{0, \bar\beta s - 1\} \,  \forall s \geq
0$;
\item[(iii)] if $x, z \in\mathbb{Z}$ have different signs and $(x, 0)
\leftrightarrow(z, s)$, then $(0, 0) \leftrightarrow(z, s)$.
\end{longlist}
\end{itproposition}

\begin{pf}
If $s\leq1$ in parts (i), (ii) or (iii), then the statements hold for
any $\bar\beta<1$, by (\ref{beta1}) and (\ref{betanew}). Hence,
from now on, we assume that $s>1$ in all three parts.
Suppose that the origin is $\beta$-expanding. Since $\Psi_{H^{(0,1)}}
\in\Gamma$, there exist sequences $\{m_n^r\}_{n \geq0}, \{m_n^l\}_{n
\geq0}$ in $\mathbb{Z}$ such that
\begin{eqnarray}\label{pcoord}
  m_0^l &=& -2,\qquad  m_0^r=2,\nonumber
  \\
|m_{n+1}^l-m_n^l| &=& |m_{n+1}^r - m_n^r| = 1 ,\nonumber
 \\[-8pt]\\[-8pt]
\Psi_{H^{(0,1)}}(m_n^l, n) &=& \Psi_{H^{(0,1)}}(m_n^r, n) = 1 ,\nonumber
\\
m_n^l &<& - \beta n < \beta n < m_n^r, \qquad n \geq0 .\nonumber
\end{eqnarray}
Define
\begin{eqnarray*}
B^l &=& \bigcup_{n=0}^\infty\bigl[(I_{m_n^l} \times KNn)\cup
J_{(m_n^l,n)}\bigr],\qquad B^r = \bigcup_{n=0}^\infty\bigl[(I_{m_n^r}
\times KNn)\cup J_{(m_n^r,n)}\bigr],
\\
B &=& B^l \cup B^r \cup(I_0 \times\{0\}).
\end{eqnarray*}

$B^l$ is a union of horizontal lines (the ``$I_m \times KNn$'''s), one
for each height level $KNn$, and rectangles of base $2M$ and height
$KN$ (the ``$J_{(m, n)}$'''s); each rectangle connects a pair of
horizontal lines. $B^l$ is thus a connected subset of $\mathbb
{R}\times[0,
+\infty)$. The same can be said about $B^r$. So, $B$ is also connected
and its complement in $\mathbb{R}\times[0, +\infty)$ has two connected
components, which will be referred to as ``above'' and ``below'' $B$.
Also, note that since $N > 2M,  \forall(x,t) \in B$, we either have
$[x-M, x] \times\{t\} \subset B$ or $[x, x+M] \times\{t\}\subset B$.
In other words, the three sets whose union defines $B$ ($B^l$, $B^r$
and $I_0 \times\{0\}$) have width larger than $M$ at any time level.

Putting together (\ref{0mscond2}), (\ref{0mscond3}), (\ref
{mscond2}), (\ref{mscond3}) and the three first conditions in (\ref
{pcoord}), we can conclude that in the trajectory $\eta^\mathbh
{1}(H^{(0,1)})$, every infected site in $(0,1)+B:= \{(z, 1+s)\dvtx (z, s)
\in B\}$ descends from $(I_{-2} \cup I_0 \cup I_2)\times\{1\}$. Then,
because of (\ref{beta2}), in the trajectory $\eta^{\mathbh{1}}(H)$,
every infected site in $(0,1) + B $ descends from $(0, 0)$.

It follows from the last condition of (\ref{pcoord}) that there exists
$0< \bar\beta<1$ such that $V(\bar\beta)$ is contained in the union
of $ (I_{-2}\cup I_0 \cup I_2) \times[0,1]$ and the area above $(0,1)+ B$.

Now, take $x$ and $z$ as in (i). Since $s > 1$ and $(z, s) \in V(\bar
\beta)$, $(z, s)$ must be above $(0, 1) + B$. So, any path starting
from $(x, 0)$ and reaching $(z, s)$ must have a point $(y, t) \in(0,
1) + B$ and thus, as we have seen, it must be the case that $(0, 0)
\leftrightarrow(y, t) \leftrightarrow(z, s)$.

Part (ii) follows from the facts that for any $s > 1, (\bar\beta s,
s)$ is to the left of $(0, 1) + B^r$, and that $\eta_s^0 \cap\{x\dvtx  (x,
s) \in B^r\}\neq\varnothing.$

Finally, take $x, z$ as in (iii) and let $\zeta$ be the path linking
$(x, 0)$ and $(z, s)$. We separately consider the two cases: there
exist $y\neq x$ and $t<1$ such that $(y, t) \in\zeta$ or not. In the
first case, (iii) follows from (\ref{beta1}). In the second case,
noting that $x$ and $z$ have different signs and $\zeta$ has
horizontal displacements of size at most $M$, and using our remarks
about $B$ being connected and its width being larger than $M$ at any
time level, we conclude that $(\gamma(t), t) \in(0,1)+B$ for some $t
\in[0, s]$. (iii) then follows from the fact that any infection in
$(0, 1) + B$ descends from $(0,0)$.
\end{pf}

\subsection{Proof of tightness of $\{\rho_t^+\}$}

Call the origin $(\beta, \gamma)$\emph{-good up to time} $T$ (resp.,
$(\beta, \gamma$)\emph{-good}) when it is both $\beta$-expanding
and $\gamma$-slow up to time $T$ (resp., $\beta$-expanding and
$\gamma$-slow). Additionally, call a point $(x, t)$  $\beta$\emph
{-expanding}, $\gamma$\emph{-slow} or $(\beta,\gamma)$\emph{-good}
when $(0, 0)$ has the corresponding property on $H^{(x, t)}$.
\begin{itlemma}\label{betagamma}
For $\gamma> 0$ sufficiently large, we have:
\begin{enumerate}[(iii)]
\item[(i)] $\mathbb{P}((0, 0) \mbox{ is } (\beta, \gamma)
\mbox
{-good}) > 0$;
\item[(ii)] $\mathbb{P}((0, 0) \mbox{ is } (\beta, \gamma)
\mbox
{-good up to time } T \mbox{ but not } (\beta, \gamma) \mbox
{-good}) \leq F \mathrm{e}^{-f T}$ for some $f, F > 0$;
\item[(iii)] given $0 \leq a < b,  \mathbb{P}((r_t, t) \mbox{
is not
} (\beta,\gamma) \mbox{-good for any } t \in[a, b]) \leq G\mathrm{e}^{-g
\sqrt{b-a}}$ for some $g,G > 0$ not depending on $a, b$.
\end{enumerate}
\end{itlemma}

\begin{pf} The only point that does not follow directly from Lemmas
\ref
{gamma} and \ref{beta} is (iii). We start proving the result when $a =
0$. Given a Harris construction $H$, define $\mu(H) = \sup\{t
\geq0\dvtx
(0, 0) \mbox{ is } (\beta, \gamma) \mbox{-good up to time } t \mbox
{ in } H\},$
\[
\sigma_1(H) = \cases{ 1 , &\quad\mbox{if } $\mu(H) < 1 $, \cr
KN(n+1) + 1 , &\quad\mbox{if } $\mu(H) \in\bigl[KNn + 1, KN(n+1) + 1\bigr) $,\cr
+\infty , &\quad\mbox{if } $\mu(H) = + \infty $,
}
\]
$\sigma_0(H) = 0 $ and $\sigma_{i+1}(H) = \sigma_i(H) + \sigma
(H^{(r_{\sigma_i(H)}, \sigma_i(H))}) \mbox{ if } i \geq1 \mbox{
and } \sigma_i(H) < +\infty.$ (The $r_{\sigma_i}$ that appears is
defined with respect to the original trajectory $\eta^{(-\infty,
0]}(H)$, with no change of coordinates.) Each $\sigma_i$ is a stopping
time for the process $t \mapsto H_t$. It follows from the strong Markov
property and translation invariance of the law of $H$ that the law of
$H^{(r_{\sigma_i}, \sigma_i)}$ conditioned to $\{\sigma_i < +\infty
\}$ is the same as that of $H$. In particular, conditioned on $\{\sigma
_i < +\infty\},$ $\sigma_{i+1} - \sigma_i$ has the law of $\sigma
_1$, which satisfies:
\begin{itemize}
\item[$\bullet$] $\mathbb{P}(\sigma_1 = +\infty) \equiv
\theta> 0$,
by (i);
\item[$\bullet$] $\mathbb{P}(T < \sigma_1 < +\infty) < \bar F
\mathrm{e}^{-\bar f T}$ for some $\bar f, \bar F > 0$, by (ii).
\end{itemize}

Let $\tau= \inf\{s\dvtx (r_s,s) \mbox{ is } (\beta,\gamma)\mbox
{-good}\}$.
Now, if $i_0$ is the first $i$ such that $\sigma_{i+1} = +\infty$, we
have $\tau\leq\sigma_{i_0}$ and
\begin{eqnarray*} \mathbb{P}(\tau> b) &\leq&\mathbb{P}(\sigma_{i_0} >
b) \leq\mathbb{P}\bigl(i_0
> \sqrt{b}\bigr) + \mathbb{P}\bigl(i_0 \leq\sqrt{b}, \sigma_{i_0} > b\bigr)
 \\
&\leq&(1-\theta)^{\sqrt{b}} + \mathbb{P}\bigl(i_0 \leq\sqrt{b}, \sigma
_{j+1} -
\sigma_j > \sqrt{b} \mbox{ for some } 1\leq j \leq i_0\bigr)
\\
&\leq&(1-\theta)^{\sqrt{b}} + \sqrt{b} \cdot\bar F \mathrm{e}^{-\bar f \sqrt
{b}} \leq G\mathrm{e}^{-g\sqrt{b}}
\end{eqnarray*}
for some suitably chosen $g, G$.

For $a>0$, repeat the proof starting from $(r_a,a)$ instead of $(0,0)$
and note that the constants $\bar f$ and $\bar F$ do not depend on $a$.
\end{pf}

\begin{pf*}{Proof of Theorem \ref{tightness} (First part)} Fix $\delta>
0$. This is the $\delta$ that takes part in our renormalization
construction, as mentioned in the paragraph after Proposition \ref
{MSprop}. We want to prove that for any~$T$, $\rho_T < L$ with
probability larger than $1-\delta$. To this end, we will proceed in
two steps. First, we will define a ``good event'' depending on $T$,
$\mathcal{H}
(T)$, with $\mathbb{P}(\mathcal{H}(T)) > 1 - \delta$. We will then
choose $L > 0$
and see that in $\mathcal{H}(T)$, every infection in $\eta^{\mathbh{1}}_T$
that is to the left of $r_T-L$ must descend from $(-\infty, 0] \times0$.

(A) \textit{Choice of the good event.}
By Lemma \ref{gamma}(i), we can choose $\gamma> 0$ such that the event
\[
\mathcal{H}_1 = \{(0, 0) \mbox{ is } \gamma\mbox{-slow}\}
\]
has probability larger than $1-\delta/3$. We can also assume that
$\gamma$ satisfies (iii) in Lemma \ref{betagamma}.

We can choose $S > 0$ such that
\[
\bigl\{\exists x \in[-S, 0] \mbox{ such that } H^{(x, 0)} \mbox{
satisfies (\ref{beta1})--(\ref
{beta2})}\bigr\}
\]
has probability larger than $1 - \delta/6$; note that this event
depends only on the Harris construction on the time interval $[0, 1]$.
Also, for any $x$, we have $\mathbb{P}(\Psi_{H^{(x, 1)}} \in\Gamma)
= \mathbb{P}
(\Psi_H \in\Gamma) > 1 - \delta/6$, by our choice of $\varepsilon$
(see the remark after Proposition \ref{MSprop}); for any $x$, this
event depends only on the Harris construction on the time interval $[1,
+\infty)$ and is thus independent of the former event. Therefore,
putting
\[
\mathcal{H}_2 = \{\mbox{there exists } x \in[-S, 0] \mbox{ such
that } (x,
0) \mbox{ is } \beta\mbox{-expanding}\},
\]
we have $\mathbb{P}(\mathcal{H}_2) > 1 - \delta/3$.

Now, choose $R > 0$ such that $\sum_{n=1}^\infty G\mathrm{e}^{-g\sqrt{R+n}} <
\delta/3$, where $g, G$ are defined in\break Lemma~\ref{betagamma}(iii).
Given $\bar R > 0$, define the time intervals
\[
I_0 = [0, \bar R],\qquad I_n = (\sup I_{n-1}, \sup I_{n-1} + R + n]\qquad \mbox{for } n \geq1,
\]
so that $I_n = (\bar R + (n-1)R + \frac{n(n-1)}{2}, \bar R + nR +
\frac{n(n+1)}{2}]$, $|I_n| = R + n$ for $n \geq1$. We now choose
$\bar R$ large enough so that\vspace*{-1pt}
%
\begin{equation}\label{overlapineq}
\forall n \geq2,   \forall t \in I_{n-1}\qquad \frac{2\bar\beta
t-S}{\bar\beta+ \gamma} > |I_{n-1} \cup I_n|.\vspace*{-1pt}
\end{equation}
Given $T > 0$, define $\bar n(T) = \sup\{n \geq1\dvtx I_n \subset[0, T]\}
$; if $I_0 \cup I_1 \nsubseteq[0, T]$, put $\bar n(T) = - \infty$.
The idea is that, given the time interval $[0, T]$, we will place the
intervals $I_n$ from top to bottom, that is, $T - I_0, T - I_1, \ldots
,$ up to the last one that fits, which will be $I_{\bar n (T)}$. Now,
define the event\vspace*{-2pt}
\begin{eqnarray*}
\mathcal{H}_3 = \mathcal{H}_3(T) &=& \{\mbox{for each } n \in[1, \bar
n(T)], \mbox{
there exists } t \in T - I_n
\\[-2pt]
&&{}\ \,\mbox{such that } (r_t, t) \mbox{ is }
(\beta,\gamma) \mbox{-good} \};\vspace*{-2pt}
\end{eqnarray*}
if $\bar n(T) = -\infty$, simply take $\mathcal{H}_3$ to be the whole space.
Now, as a consequence of Lemma \ref{betagamma}(iii), we obtain\vspace*{-2pt}
\begin{eqnarray*}
\mathbb{P}(\mathcal{H}_3(T)) &\geq&1 - \sum
_{n=1}^{\bar n(T)}\mathbb{P}\bigl((t,
r_t) \mbox{ is never }(\beta,\gamma) \mbox{-good when } t \in T -
I_n\bigr)
\\[-2pt]
&\geq&1 - \sum_{n=1}^\infty G\mathrm{e}^{-g\sqrt{|I_n|}} = 1 - \sum
_{n=1}^\infty G\mathrm{e}^{-g\sqrt{R+n}} > 1 - \delta/3.
\end{eqnarray*}

In conclusion, if $\mathcal{H}(T) = \mathcal{H}_1 \cap\mathcal{H}_2
\cap\mathcal{H}_3(T)$, then
$\mathbb{P}(\mathcal{H}(T)) > 1 - \delta$ for any $T$.

(B) \textit{Choice of} $L$ \textit{and proof that the interface area
is smaller than} $L$ \textit{in the good event}. Let $L = \gamma(R +
\bar R + 1) + S$; note that $L$ does not depend on $T$. We first treat
the case $T \leq\bar R + R + 1$. We might omit it: since $\sup_{t\leq
T}|\rho_t| < \infty$ almost surely, it suffices to prove its
tightness in $[T, +\infty)$ for sufficiently large $T$. However, we
find that this case illustrates the main idea of the proof without the
technical complications that appear in the general picture.

Let $V = V(\bar\beta) = \{(z, s) \in\mathbb{Z}\times[0, +\infty)\dvtx
-\bar
\beta s \leq z \leq\bar\beta s\}$, where $\bar\beta>0 $ is such
that the conclusion of part (i) of Proposition \ref{barrier} holds.
Given $A \subset\mathbb{Z}\times[0, +\infty)$ and $t \geq0$,
define $\Pi
_t(A) = \{z\dvtx (z, t) \in A\}$.

Fix $H \in\mathcal{H}(T)$. Since $H \in\mathcal{H}_2$, we can take
$x \in[-S, 0]$
such that $(x, 0)$ is $\beta$-expanding. Also, since $H \in\mathcal{H}_1,
(0,0)$ is $\gamma$-slow and, in particular, $r_T < \gamma T$. Thus,
\begin{eqnarray*} r_T - L &<& \gamma T - L \leq\gamma(\bar R + R + 1)
- \gamma(\bar R + R + 1 ) - S
\\
&=& - S < x + \bar\beta T < \sup\Pi
_T\bigl((x,0) + V\bigr)+1;
\end{eqnarray*}
the $+1$ is required because $x + \bar\beta T$ may not be an integer.
Assume that for $y > 0$ and $w$ satisfying $r_T - w > L$, we have $(y,
0) \leftrightarrow(w, T)$. Note that $w < r_T - L \leq\sup\Pi_T((x,
0) + V)$. If $w \in\Pi_T((x, 0) + V)$, then it follows from
Proposition \ref{barrier}(i) and translation invariance that $(x, 0)
\leftrightarrow(w, T)$. If $w < \inf\Pi_T((x, 0) + V)$, then $w$ and
$y$ are in opposite sides of $x$ and it follows from Proposition \ref
{barrier}(ii) and translation invariance that $(x, 0) \leftrightarrow
(w, T)$. This shows that any infection in $(-\infty, r_T - L] \times
T$ that descends from $[1, +\infty) \times0$ must also descend from
$(-\infty, 0]\times0$, completing the proof of this case.

Before starting the other case, we make some trivial remarks. Suppose
$(a, s), (b, t) \in\mathbb{Z}\times[0, +\infty)$ are such that $a
\leq b$
and $s < t$. Let $\zeta^*$ be the smallest value of $\zeta$ at which
$\Pi_\zeta((a, s) + V) \cap\Pi_\zeta((b, t) + V) \neq\varnothing$.
$\zeta^*$ is either $t$ (in the case $(b,t) \in(a, s)+V$) or the time
of intersection of the lines $\zeta\mapsto a +\bar\beta(\zeta-s)$
and $\zeta\mapsto b - \bar\beta(\zeta- t)$, that is,
%
\begin{equation}\label{zetastar}\zeta^*((a,s),(b,t)) = \max\biggl\{
t, \frac{b-a+\bar
\beta(t+s)}{2\bar\beta}\biggr\}.
\end{equation}
Also,
%
\begin{equation}\label{overlap} \zeta> \zeta^*((a,s),(b,t))
\quad\Longrightarrow\quad \Pi
_\zeta\bigl((a,s)+V\bigr) \cup\Pi_\zeta\bigl((b,t)+V\bigr) \mbox{ is an interval}.
\end{equation}

Now, take $T > \bar R + R + 1$ and $H \in\mathcal{H}(T)$. Again, $(0,
0)$ is
$\gamma$-slow and there exists $x \in[-S, 0]$ such that $(x, 0)$ is
$\beta$-expanding. Also, since $H \in\mathcal{H}_3(T)$, there exist
$t_1 \in
T-I_{\bar n}, t_2 \in T - I_{\bar n -1}, \ldots, t_{\bar n} \in T -
I_1$ such that $(r_{t_i}, t_i)$ is $(\beta,\gamma)$-good for $i = 1,
\ldots, \bar n$. Note that since $(0 , 0)$ and each $(r_{t_i}, t_i)$
is $\gamma$-slow, we have
%
\begin{eqnarray}\label{beforegamma}
 r_{t_1} &\leq&\gamma t_1 ,\nonumber
 \\
 r_{t_{n+1}} &\leq& r_{t_n}
+ \gamma(t_{n+1}-t_n),\qquad n = 1, \ldots, \bar n - 1 ,
\\
r_T &\leq&
r_{t_{\bar n}} + \gamma(T - t_{\bar n}) ,\nonumber
\end{eqnarray}
and by Proposition \ref{barrier}(ii) and translation invariance, we have
%
\begin{eqnarray}\label{afterbeta}
 r_{t_1} &\geq& x,\nonumber
 \\[-8pt]\\[-8pt]
 r_{t_{n+1}} &\geq& r_{t_n},\qquad n = 1,
\ldots, \bar n - 1.\nonumber
\end{eqnarray}

We claim that the cones $(r_{t_i}, t_i) + V$ each overlap with their
neighbors before time $T$, that is,
%
\begin{eqnarray}\label{overneighbor}
\zeta^*((x, 0), (r_{t_1}, t_1)) &<& T ,\nonumber
\\[-8pt]\\[-8pt]
\zeta^*((r_{t_i},
t_i),(r_{t_{i+1}}, t_{i+1})) &<& T, \qquad i = 1, \ldots, \bar n - 1 .\nonumber
\end{eqnarray}
Let us prove the first expression in (\ref{overneighbor}). If
$(r_{t_1}, t_1) \in(x, 0) + V$, then $\zeta^*((x,0), (r_{t_1}, t_1))
= t_1 < T$. Assume that $(r_{t_1}, t_1) \notin(x, 0)+V$. Since $-S
\leq x < r_{t_1} \leq\gamma t_1$, we have $r_{t_1}-x \leq\gamma t_1 +
S$. Also,
\[
0 \in T - I_{\bar n + 1}\quad \Longrightarrow\quad T \in I_{\bar n +1} \quad\stackrel
{(\ref{overlapineq})}{\Longrightarrow}\quad \frac{2 \bar\beta T -
S}{\bar\beta+ \gamma} > |I_{\bar n +1}\cup I_{\bar n + 2}|>|I_{\bar
n}\cup I_{\bar n + 1}|
\]
and since we also have that $t_1 \in T - I_{\bar n}$, we obtain $t_1 =
t_1 - 0 < |I_{\bar n} \cup I_{\bar n +1}| < \frac{2 \bar\beta T -
S}{\bar\beta+ \gamma}$. Putting these inequalities together and
using (\ref{zetastar}), we get
\begin{eqnarray*}
\zeta^*((x, 0), (r_{t_1}, t_1)) &=& \frac{r_{t_1}-x + \bar\beta
t_1}{2\bar\beta} \leq\frac{\gamma t_1 + S + \bar\beta t_1}{2 \bar
\beta}
\\
&<& \frac{S}{2\bar\beta} + \frac{2 \bar\beta T - S}{\bar
\beta+ \gamma} \cdot\frac{\gamma+ \bar\beta}{2 \bar\beta} = T.
\end{eqnarray*}
For the second expression in (\ref{overneighbor}), if $(r_{t_{i+1}},
t_{i+1}) \in(r_{t_i}, t_i) + V$, then $\zeta^*((r_{t_i}, t_i),
(r_{t_{i+1}}, t_{i+1})) = t_{i+1} < T$. Assume that $(r_{t_{i+1}},
t_{i+1}) \notin(r_{t_i}, t_i) + V$ and write
\begin{eqnarray*}
\zeta^*((r_{t_i}, t_i),(r_{t_{i+1}}, t_{i+1}))&=&\frac
{r_{t_{i+1}}-r_{t_i}+ \bar\beta(t_{i+1}+t_i)}{2\bar\beta}\leq
\frac{\gamma(t_{i+1}-t_i)+\bar\beta(t_{i+1}+t_i)}{2\bar\beta}
\\[2pt]
&=&\frac{\gamma(t_{i+1}-t_i)}{2\bar\beta}+\frac{t_{i+1}+t_i}{2} .
\end{eqnarray*}
Since $t_i \in T-I_{\bar n - i + 1}$ and $t_{i+1}\in T-I_{\bar n - i}$,
we have $t_{i+1}-t_i \leq\vert I_{\bar n-i}\cup I_{\bar n-i+1} \vert
\leq\frac{2 \bar\beta t - S}{\bar\beta+ \gamma} $
for any $t\in I_{\bar n-i}$, by (\ref{overlapineq}). In particular,
this holds for $t= T - t_{i+1}$.
Therefore,
\begin{eqnarray*} \zeta^*((r_{t_i}, t_i),(r_{t_{i+1}}, t_{i+1}))&\leq&
\frac{\gamma}{2 \bar\beta} \cdot\frac{2\bar\beta(T-
t_{i+1})-S}{\bar\beta+ \gamma}
\leq\frac{2\bar\beta
(T-t_{i+1})-S}{2\bar\beta} +\frac{t_{i+1}+t_i}{2}
\\[2pt]
&\leq& T-t_{i+1}+
\frac{t_{i+1}+t_i}{2}\leq T.
\end{eqnarray*}

Now, define the union of cones $U = [\bigcup_{n=1}^{\bar n}
((r_{t_n}, t_n)+V) ] \cup[(x, 0) + V]$. Using (\ref{overlap})
and (\ref{overneighbor}), we conclude that $\Pi_T(U)$ is an interval.

Since $t_{\bar n} \in T - I_1$, we have $T - t_{\bar n} \leq\bar R + R
+ 1$. Also, using the last inequality in (\ref{beforegamma}), we obtain
\[
r_T < r_{t_{\bar n}} + \gamma(T - t_{\bar n}) < r_{t_{\bar n}} +
\gamma(\bar R + R + 1),
\]
so, using $L = \gamma(\bar R + R + 1) + S$, we have
%
\begin{equation}\label{ceilingspace}
r_T - L < r_{t_{\bar n}} + \gamma(\bar R + R +1) - \gamma(\bar R + R
+1) - S < r_{t_{\bar n}} < \sup\Pi_T(U).
\end{equation}

As before, take $y > 0$ and $w$ satisfying $r_T-w > L$ and $(y,0)
\leftrightarrow(w, T)$. Since $w < r_T - L < \sup\Pi_T(U)$ and $\Pi
_T(U)$ is an interval, there are two possibilities:
\begin{enumerate}[(a)]
\item[(a)] $w \in\Pi_T(U)$\\
In this case, by the definition of $U$, we either have $w \in\Pi
_T((x,0) + V)$ (hence $(x, 0) \leftrightarrow(w, T)$, as we already
saw) or $w \in\Pi_T((r_{t_i}, t_i) + V)$ for some $i$.
In this last case, there exists $z$ such that $(z,t_i)\leftrightarrow
(w,T)$ and hence, by part (i) of Proposition (\ref{barrier}),
$(r_{t_i}, t_i) \leftrightarrow(w, T)$, which implies that $(-\infty,
0] \times\{0\} \leftrightarrow(w, T)$.

\item[(b)] $w < \inf\Pi_T(U)$\\
By the same argument that was used in the case $T < \bar R + R + 1$, we
have $(x, 0) \leftrightarrow(w, T)$.
\end{enumerate}

In conclusion, in any case, we have $(-\infty, 0] \times0
\leftrightarrow(w, T)$. We have thus shown that any point $(w, T)$
that is connected to $[1, +\infty)$ but not to $(-\infty, 0]$ must be
to the right of $(r_T -L, T)$, that is, that $\rho_T < L$, as required.
\end{pf*}

\section{Tightness of $\{\rho_t^-\}$}
\label{S3}

In the following lemma, we will reuse the renormalization structure
built in the last section. We fix an arbitrary $\beta\in(0, 1)$ and
$k, K$ as in Proposition \ref{MSprop}, then choose a closure density
$\varepsilon$ such that the event $\Gamma$ of Lemma \ref{betap} has
positive probability. Finally, we choose $N$ such that $\Psi_H$ has
closure density below $\varepsilon$ (again as in Proposition \ref{MSprop}).

\begin{itlemma}\label{maximum} For any $\sigma> 0$, there exists $L >
0$ such that for
any $T > 0$,
%
\begin{equation}\label{boundmax} \mathbb{P}(\mbox{there exists } t
\leq T \mbox{ such that }
r_t > r_T + L) < \sigma.
\end{equation}
\end{itlemma}

\begin{pf}
As in the proof of Theorem \ref{tightness}, we will define an event
$\mathcal{G}= \mathcal{G}_1 \cap\mathcal{G}_2 \cap\mathcal
{G}_3(T)$ such that $\mathbb{P}(\mathcal{G}) >
1-\sigma$ and choose an appropriate $L > 0$; we will then show that in
$\mathcal{G}$, we have
%
\begin{equation}\label{Lbound} r_t < r_T + L  \qquad \forall t \leq T.
\end{equation}

The first event is the same as before: $\mathcal{G}_1 = \{(0,0) \mbox
{ is }
\gamma\mbox{-slow}\}$, with $\gamma$ chosen so that this has
probability $> 1 - \sigma/3$ (see Lemma \ref{gamma}). Put $\mathcal
{G}_2 = \{
r_t > - S  \, \forall t \geq0\}$ with $S > 0$ chosen such that this has
probability greater than $ 1 - \sigma/ 3$; this is possible because
$\inf\{r_t \dvtx t \geq0\} > -\infty$ almost surely.

Increasing $\gamma$ so that the conclusions of Lemma \ref{betagamma}
hold, we may choose $R > 0$ such that $\sum_{n=0}^\infty G \mathrm{e}^{-g\sqrt
{R+n}} < \sigma/3$, where $g, G$ are as in part (iii) of Lemma \ref
{betagamma}. We then put
\[
I_0 = [0, R),\qquad  I_n = [\sup I_{n-1}, \sup I_{n-1} + R +n)\qquad \mbox{for } n
\geq1,
\]
so that $I_n = [nR + \frac{(n-1)n}{2},(n +1)R + \frac{(n+1)n}{2})$,
$|I_n| = R + n$ when $n \geq0$. We also put $\bar n(T) = \sup\{n \geq
0\dvtx I_n \subset[0, T]\}$; if $I_0 \nsubseteq[0, T]$, then put $\bar
n(T) = - \infty$. Next, define
\[
\mathcal{G}_3(T) = \{\mbox{for each } n \in[0, \bar n(T)], \mbox{ there
exists } t \in T - I_n \mbox{ such that } (r_t, t) \mbox{ is } (\beta
, \gamma)\mbox{-good}\};
\]
if $\bar n(T) = -\infty$, take $\mathcal{G}_3$ to be the whole space.
By the
choice of $R$ and Lemma \ref{betagamma}(iii), $\mathbb{P}(\mathcal
{G}_3(T)) > 1 -
\sigma/3$. Thus, $\mathbb{P}(\mathcal{G}) > 1 - \sigma$, as required.

Let us recall that
%
\begin{equation}\label{goodbounds} (r_s, s) \mbox{ is } (\beta,
\gamma)\mbox{-good},\qquad
s' > s\quad \Longrightarrow\quad r_s + \bar\beta(s'-s) -1 \leq r_{s'} \leq r_s
+ \gamma(s'-s) ,
\end{equation}
where $\bar\beta$ is defined in Proposition \ref{barrier}. Choose
$L$ such that
%
\begin{eqnarray}\label{L1}L &\geq&\gamma(2R + 1) + S \quad \mbox{and}
\\
\label{L2} L &\geq& S + \gamma(2R + 2n +1)-\bar\beta
\biggl(nR+\frac
{n(n-1)}{2}\biggr) + 1  \qquad \forall n \geq0.
\end{eqnarray}

We proceed to prove that (\ref{Lbound}) is satisfied in $\mathcal
{G}$. Fix $0
< t < T$. We deal with three cases:
\begin{itemize}[$\bullet$]
\item[$\bullet$] $t < T \leq2R + 1$. Since the origin is
$\gamma$-slow, we have $r_t \leq\gamma t \leq\gamma(2R + 1).$ Since
we are in $\mathcal{G}_2$, we have $r_T > -S$. Therefore, $r_T + L >
-S + L
\stackrel{(\ref{L1})}{\geq} \gamma(2R + 1) \geq r_t$.
\item[$\bullet$] $T > 2R + 1, t \in(T-I_{\bar n})\cup
(T-I_{\bar n + 1})$ (the point being that $t$ is close to zero, so
there does not necessarily exist a $(\beta, \gamma)$-good point below
$(r_t, t)$). By the definition of $\bar n$, we have $0 \in I_{\bar n +
1}$, so $t < |I_{\bar n} \cup I_{\bar n + 1}| = 2R + 2 \bar n + 1$ and
$r_t < \gamma t < \gamma(2R + 2 \bar n + 1)$. Also, by the definition
of $\mathcal{G}_3$, there exists $t^* \in T - I_{\bar n}$ such that $(r_{t^*},
t^*)$ is $(\beta, \gamma)$-good.\vspace*{1pt} $t^* \in T - I_{\bar n}$ implies
that $T- t^* \geq\inf I_{\bar n} = \bar n R + \frac{(\bar n - 1)\bar
n}{2}$. We then have
\begin{eqnarray*}
r_T + L &\stackrel{(\ref{goodbounds})}{\geq}& r_{t^*} + \bar\beta
(T-t^*) + L - 1 > -S + \bar\beta\biggl(\bar n R + \frac{(\bar n -
1)\bar n}{2}\biggr) + L - 1
\\
&\stackrel{(\ref{L2})}{\geq}& \gamma(2R +
2 \bar n + 1) \geq r_t.
\end{eqnarray*}
\item[$\bullet$] $T > 2R + 1$, $t \in T - I_n$ with $n < \bar
n$. Here, $n + 1 \leq\bar n$, so there exists $t^* \in T - I_{n+1}$
such that $(r_{t^*}, t^*)$ is $(\beta, \gamma)$-good. Note that $t >
t^*, t- t^* < |I_n \cup I_{n+1}| = 2R + 2n + 1$, so (\ref{goodbounds}) gives
%
\begin{equation}\label{rtbound}r_t \leq r_{t^*} + \gamma(t-t^*) \leq
r_{t^*} + \gamma
(2R + 2n +1).
\end{equation}
On the other hand, $T - t^* \geq|I_0 \cup\cdots\cup I_n| = (n+1)R +
\frac{(n+1)n}{2}$, so
\begin{eqnarray*}r_T + L &\stackrel{(\ref{goodbounds})}{\geq}& r_{t^*} +
\bar\beta(T - t^*) + L - 1 \geq r_{t^*} + \bar\beta\biggl(nR +
\frac{(n+1)n}{2}\biggr) + L - 1
\\[-2pt]
&\stackrel{(\ref{L2})}{\geq}&
r_{t^*} + \gamma(2R + 2n + 1) \stackrel{(\ref{rtbound})}{\geq} r_t.
\end{eqnarray*}
\end{itemize}
\upqed\end{pf}

\begin{itlemma}\label{rightpath} For any $\sigma> 0$, there exists $L
> 0$ such that
for any $T > 0$,
%
\begin{equation}\label{pathexists}
\mathbb{P}\bigl([0, +\infty) \times0 \leftrightarrow[0, L] \times T
\mbox{
inside } (0, +\infty)\bigr) > 1-\sigma.
\end{equation}
\end{itlemma}

This follows from the fact that $r_t$ has positive asymptotic speed and
a simple duality argument; we omit the proof.

For $T > 0$, define $q_T = \max\{r_t\dvtx 0 \leq t \leq T\}$. We now
proceed to complete the proof of Theorem \ref{tightness}.

\begin{pf*}{Proof of Theorem \ref{tightness} (Second part)} Fix $\delta
> 0$. By Lemmas \ref{maximum} and \ref{rightpath}, we can obtain
$L_1, L_2 > 0$ such that
\begin{eqnarray*}
\mathbb{P}(q_T \leq r_T + L_1) &>& \sqrt{1 - \delta},
\\
\mathbb{P}\bigl([0, +\infty) \times0 \leftrightarrow[0, L_2] \times T
\mbox{
inside } \{(x, t)\dvtx x \geq0\}\bigr) &>& \sqrt{1 - \delta}.
\end{eqnarray*}
Put $L = L_1 + L_2 + M$. For any $T>0$, we have
\begin{eqnarray*}
\mathbb{P}(\rho_t \geq-L)
&=&\mathbb{P}\bigl((0, +\infty) \times0 \leftrightarrow(r_T, r_T+L]\times
T\bigr)
\\
 &\geq&
\mathbb{P}\bigl(q_T \leq r_T + L_1,
[r_T+L_1+M+1, +\infty)\times0
\\
&&{}\quad\hspace*{-3pt} \leftrightarrow[r_T + L_1 + M +1, r_T+L] \times T  \mbox{ inside }
[r_T + L_1 + M +1, + \infty)
\bigr)
\\
&=& \sum_{x=-L_1}^{+\infty} \mathbb{P}\bigl(r_T=x, q_T \leq x+ L_1,
[x+L_1+M+1, +\infty)\times0
\\
&&{}\qquad\qquad\hspace*{-7pt}\leftrightarrow[x+L_1+M+1, x+L]\times T
 \mbox{ inside } [x + L_1 + M +1, + \infty)
\bigr).
\end{eqnarray*}
(The sum starts at $-L_1$ because $q_T \geq0$, so we can only have
$q_T \leq r_T+L_1$ when $r_T \geq-L_1$.) Now, in each of the above
probabilities, the first two events depend on the Harris construction
on the set $(-\infty, x+L_1+M] \times[0, +\infty)$, whereas the
third event depends on the Harris construction on $[x+L_1+M+1, +\infty
) \times[0, +\infty)$. They are thus independent and the sum becomes
\begin{eqnarray*}
&&\sum_{x=-L_1}^{+\infty}\mathbb{P}(r_T = x, q_T \leq x+L_1)
\cdot\mathbb{P}\bigl( [x+L_1+M+1, +\infty)\times0
\\[-10pt]
&&\quad\hspace*{131pt}\leftrightarrow[x+L_1+M+1, x+L]\times T
\\
&&\quad\hspace*{131pt} \mbox{inside } [x + L_1 +
M +1, + \infty)
\bigr)
\\
&&\quad=\mathbb{P}\bigl([0, +\infty) \times0 \leftrightarrow[0, L_2]\times T
\mbox{
inside } [0, +\infty)\bigr)\cdot\sum_{x=-L_1}^{+\infty} \mathbb{P}(r_T
= x, q_T
\leq x+L_1)
 \\
&&\quad =\mathbb{P}\bigl([0, +\infty)\times0 \leftrightarrow[0, L_2]\times T
\mbox{
inside } [0, +\infty)\bigr)\cdot\mathbb{P}(q_T \leq r_T +L_1) > 1 -
\delta,
\end{eqnarray*}
completing the proof.
\end{pf*}


\printhistory

\end{document}